\documentclass[11pt,a4paper]{amsarticle}
\usepackage[english]{babel}
\usepackage[latin1]{inputenc}
\usepackage{amsmath,latexsym,amssymb}
\usepackage{amsfonts}
\usepackage{amsthm}
\usepackage{url}
\usepackage{color}

\newtheorem{thm}{Theorem}[section]
\newtheorem{cor}[thm]{Corollary}
\newtheorem{defi}[thm]{Definition}

\newtheorem{prop}[thm]{Proposition}
\newtheorem{nota}[thm]{Notation}

\theoremstyle{remark}
\newtheorem{rem}[thm]{\bf Remark}


\begin{document}

\title[Infinite-dimensional manifolds as ringed spaces]{Infinite-dimensional manifolds \\ as ringed spaces}

\author{Michel Egeileh} 
\address{%
Facult\'e des sciences\\
Universit\'e Saint-Joseph\\
B.P. 11-514 Riad El-Solh\\ 
Beyrouth 1107 2050\\
\mbox{Liban}
}
\email{michel.egeileh@usj.edu.lb}

\author{Tilmann Wurzbacher} 

\address{%
Institut \'Elie Cartan Lorraine\\
Universit\'e de Lorraine et
C.N.R.S.\\
57045 Metz, France}

\email{tilmann.wurzbacher@univ-lorraine.fr}

%

\begin{abstract}
We analyze the possibility of defining infinite-dimensional manifolds as ringed spaces. More precisely, we consider three definitions of manifolds modeled on locally convex spaces: in terms of charts and atlases, in terms of ringed spaces, and in terms of functored spaces, as introduced by  Douady in
his thesis. It is shown that for large classes of locally convex model spaces (containing Fr\'echet spaces and duals of Fr\'echet-Schwartz spaces), the three definitions are actually equivalent. The equivalence of the definition via charts with the definition via ringed spaces is based on the fact that for the classes of model spaces under consideration, smoothness of maps turns out to be equivalent to their {\it scalarwise smoothness} (that is, the smoothness of their composition with smooth real-valued functions). \\
\end{abstract}

\maketitle



\section{Introduction}

Finite-dimensional manifolds, whether smooth, real- or complex-analytic, are commonly defined via charts and atlases. The other standard way of defining them relies on a dual point of view, focussing on the functions rather than on the points themselves, and this is achieved via a sheaf-theoretical approach. More precisely, a smooth $n$-dimensional manifold $M$ is then defined as a locally ringed space $(M_0,\mathcal{O}_M)$ that is locally isomorphic to the locally ringed space $(\mathbb{R}^n,\mathcal{C}^{\infty}_{\mathbb{R}^n})$.
The sheaf-theoretical approach is hardly avoidable when one wants to deal with singular generalizations of manifolds (varieties or schemes for instance), or with ``non-reduced situations", such as supermanifolds, where the rings of ``functions" have nilpotents. In this last example, a section of the structural sheaf is not determined by its values on the points of the underlying topological space, which makes the sheaf-theoretical approach particularly relevant in defining supermanifolds.  \\

In finite dimensions, the two definitions of manifolds (via atlases and as certain locally ringed spaces respectively) are well-known to be equivalent. In infinite dimensions, the situation is quite different. Infinite-dimensional manifolds, whether locally modeled on Banach spaces, Fr\'echet spaces or general locally convex spaces, have almost always been defined in terms of charts and atlases. One reason for that is the belief, following the thesis of Douady \cite{dou}, that the sheaf of scalar-valued functions does not give sufficient information to define the morphisms (contrary to the finite-dimensional case, where defining the smooth functions valued in $\mathbb{R}$ suffices to determine the morphisms valued in $\mathbb{R}^k$ for every natural number $k$). In \cite{maz}, Mazet defines a category of infinite-dimensional analytic spaces, precisely in terms of ringed spaces. However, his category leads to pathologies (such as the sum of two analytic maps not necessarily being analytic). Douady avoids these pathologies by introducing a third approach for capturing the notion of space, which he uses to define his category of Banach analytic spaces. Namely, given a category $\mathcal{C}$, Douady defines a {\it $\mathcal{C}$-functored space} $X$ to be a pair $(X_0,\mathcal{O}^{\mathcal{C}}_X)$ where $X_0$ is a topological space, and $\mathcal{O}^{\mathcal{C}}_X$ is a covariant functor from $\mathcal{C}$ to the category of sheaves of sets on $X_0$. In this way, for every object $F$ in $\mathcal{C}$ (thought of as a possible target), one associates a sheaf of sets $\mathcal{O}^{\mathcal{C}}_X(F)$ (thought of as the sheaf of $F$-valued morphisms on $X_0$). Compared to a ringed space, a functored space encodes already in its ``structural functor" the definition of the morphisms valued in any target space (from a certain category). \\

The functored space approach obviously adds a supplementary ``technical layer", which can be felt already when defining the local models for Banach analytic spaces (as functored spaces). Thus, unless the recourse to functored spaces is absolutely necessary, it is preferable to deal with the more traditional setting of ringed spaces, i.e. to associate only a single structure sheaf instead of a
sheaf-valued functor to each space. We are thus lead to the question whether the insufficiency of the sheaf of scalar-valued functions pointed out by Douady (and the related pathologies) appears also in the non-singular setting of infinite-dimensional smooth manifolds.\\

In this paper, we address this question by observing that the obstruction to define infinite-dimensional manifolds as ringed spaces boils down to the failure of a {\it scalarwise smooth} map between open sets of locally convex spaces to be smooth. More precisely, given two locally convex spaces $E$ and $F$, and given an open subset $U$ of $E$, we use a standard notion of smoothness for maps $\Phi:U\longrightarrow F$ ({\it cf.} Definition \ref{BHMN}), going back at least to Bastiani and adopted notably by Hamilton, Milnor and Neeb. Then, a map $\Phi$ is said to be {\it scalarwise smooth} if the function $f\circ\Phi:U\longrightarrow\mathbb{R}$ is smooth for every smooth function $f:F\longrightarrow \mathbb{R}$. The chain rule implies clearly that smooth maps are scalarwise smooth. The converse is easily seen to be true in finite dimensions (just take the linear forms $e_i^*$ dual to a basis $\{e_i \;;\,  1\leq i \leq \hbox{dim}\,F\}$). In infinite dimensions, the converse is non-trivial. Our first main result is to prove it for large classes of locally convex spaces ({\it cf.}  Theorem \ref{scalarwise smoothness implies smoothness} in the body of the article.)\\

\noindent {\bf Theorem A.}
{\it Let $E$ and $F$ be locally convex space spaces, and $U$ an open subset of $E$. Assume that for every $n\geq 1$, the $\hbox{\upshape{c}}^{\infty}$-topology on $E^n$ is the same as the product topology for the given topology on $E$ and that F is Mackey-complete. Let $\Phi:U\longrightarrow F$ be a continuous map. Then the following are equivalent:
\begin{enumerate}
\item [(i)] $\Phi$ is smooth.
\item [(ii)] For every open subset $V$ of $F$ containing $\Phi(U)$ and every $f\in\mathcal{C}^{\infty}(V,\mathbb{R})$, we have $f\circ\Phi\in \mathcal{C}^{\infty}(U,\mathbb{R})$.

\item [(iii)] For every $f\in \mathcal{C}^{\infty}(F,\mathbb{R})$, we have $f\circ\Phi\in \mathcal{C}^{\infty}(U,\mathbb{R})$ (i.e. $\Phi$ is scalarwise smooth).
\item [(iv)] For every $\ell\in F'$, we have $\ell\circ\Phi\in \mathcal{C}^{\infty}(U,\mathbb{R})$ (i.e. $\Phi$ is weakly smooth). 
\end{enumerate} }

\vskip0.1cm
Note that the assumption on $E$ is satisfied, e.g., by Fr\'echet spaces and duals of Fr\'echet-Schwartz spaces, and that every complete locally convex space is Mackey-complete. Note furthermore that the $\hbox{\upshape{c}}^{\infty}$-topology put forward by Kriegl and Michor in the fundamental work \cite{km} is also called the Mackey-closure topology.  \\

As underlined by the preceding result, in infinite dimensions one often has to single out a class $\mathcal E$ of locally 
convex model spaces. 
A Hausdorff space $M_0$, together with a smooth atlas $\mathcal A$ of charts taking values in spaces of this class, will be called a
{\it smooth $\mathcal E$-manifold}, whereas a {\it structure sheaf-smooth $\mathcal E$-manifold} is a locally ringed space 
$(M_0,{\mathcal O}_M)$ that is locally isomorphic as such to open sets of spaces in the model class (together with their natural sheaves
of smooth scalar-valued functions). In this language our next main result reads as follows 
({\it cf.}  Theorem  \ref{smooth equals ringed morphism} and 
Corollary \ref{manifolds = structure sheaf-smooth manifolds}).\\

\noindent {\bf Theorem B.}
{\it Let $\mathcal{E}$ be the class of Mackey-complete locally convex spaces $E$ such that for every $n\geq 1$, the $\hbox{\upshape{c}}^{\infty}$-topology on $E^n$ is the same as the product topology for the given topology on $E$, and let $\mathcal{M}$ be the class of Mackey-complete locally convex spaces. Then:
\begin{enumerate}
\item [(i)] If $M=(M_0,\mathcal{A})$ is a smooth $\mathcal{E}$-manifold, $N=(N_0,\mathcal{B})$ a smooth $\mathcal{M}$-manifold, and $\Phi:M_0\longrightarrow N_0$ a continuous map, then $\Phi$ is smooth if and only if $\Phi:(M_0,\mathcal{C}^{\infty}_M)\longrightarrow (N_0,\mathcal{C}^{\infty}_N)$ is a morphism of locally ringed spaces.
\item [(ii)] For every structure sheaf-smooth $\mathcal{E}$-manifold $(M_0,\mathcal{O}_M)$, there is a canonical maximal atlas $\mathcal{A}$ on $M_0$ such that $M=(M_0,\mathcal{A})$ is a smooth $\mathcal{E}$-manifold fulfilling $\mathcal{C}^{\infty}_M=\mathcal{O}_M$. Furthermore, the maximal atlas $\mathcal{A}$ is uniquely determined by the condition $\mathcal{C}^{\infty}_M=\mathcal{O}_M$.
\end{enumerate} }

\vskip0.1cm

We complete our comparison of the various definitions of infinite-dimen\-sional manifolds by showing that the definition via charts and atlases is equivalent to the one based on functored spaces.\\

\noindent Furthermore, we prove that for certain infinite-dimensional manifolds, smoothness of a continuous map is characterized in terms of pulling back globally defined smooth functions to globally defined smooth functions ({\it cf.}  Theorem  \ref{local = global}).
\vskip0.15cm

\noindent {\bf Theorem C.}
{\it Let  $\mathcal{E}$ be the class of locally convex spaces $E$ such that for every $n\geq 1$, the $\hbox{\upshape{c}}^{\infty}$-topology on $E^n$ is the same as the product topology for the given topology on $E$. Also, let $M=(M_0,\mathcal{A})$ a smooth $\mathcal{E}$-manifold, and $N=(N_0,\mathcal{B})$ a smooth regular manifold modeled on a 
nuclear Fr\'echet space $E$. Finally, let $\Phi:M_0\longrightarrow N_0$ be a continuous map. Then the following are equivalent: 
\begin{enumerate}
\item [(i)] $\Phi$ is smooth.
\item [(ii)] For every $f\in\mathcal{C}^{\infty}_N(N_0)$, we have $f\circ\Phi\in\mathcal{C}^{\infty}_M(M_0)$.
\end{enumerate}  }

\vskip0.1cm
We conclude by ascertaining that this notably holds true if the target manifold $N$ is the space of smooth maps between finite-dimensional 
manifolds (with compact source manifold) and the manifold $M$ is modeled on Fr\'echet spaces. (In this article, all ``compact manifolds"
are closed.)\\

Our paper is organized as follows. In Section 2, we present a proof of the special case of Theorem \ref{scalarwise smoothness implies smoothness}, when the domain is $\mathbb{R}$. More precisely, we show that for a Mackey-complete locally convex space $E$, a curve $c:\mathbb{R}\longrightarrow E$ is smooth if and only if it is scalarwise smooth. While this result and its idea of proof are not new (compare \cite{km}), they are crucial to our proof of Theorem \ref{scalarwise smoothness implies smoothness}: we recall them in a concise but self-contained way for the convenience of the reader, which gives us also the opportunity to introduce our notations for the rest  of the paper. In Section 3, we recall the definition of smooth maps that we will be using, and prepare and prove Theorem \ref{scalarwise smoothness implies smoothness}, using as an intermediate step the calculus of {\it convenient smoothness} studied by Kriegl and Michor. In Section 4, we present in detail the three definitions of infinite-dimensional manifolds under investigation here, and prove the comparison results (Theorem \ref{smooth equals ringed morphism} and Corollary \ref{manifolds = structure sheaf-smooth manifolds}) mentioned above. In Section 5, we prove the global characterization of smoothness mentioned above, and discuss the important example class of mapping spaces with compact source. \\
 
\noindent \emph{Acknowledgements.} This research was supported by the Ruhr-Universit\"at Bochum and the SFB/TR 12 of the Deutsche Forschungsgemeinschaft.
We would like to thank Alexander Alldridge for discussions on Douady's ``espaces fonct\'es". We would like to extend our thanks to the referee for her/his extremely useful remarks.

\section{Smoothness of curves}

The goal of this section is to recall the proof of the following result: if $c$ is a map from $\mathbb{R}$ to a complete real (Hausdorff) locally convex space $E$, and $E'$ is the continuous dual of $E$, then smoothness of $\ell\circ c$ for every $\ell\in E'$ implies the smoothness of $c$ (Theorem \ref{weak smoothness of  curves}). While this is obvious if $E$ is finite-dimensional (it is enough to take the projections $e_i^*:E\longrightarrow\mathbb{R}$ where $\{e_i \;;\,  1\leq i \leq \hbox{dim}\,E\}$ is an arbitrary basis of $E$), proving that it remains true in the general case requires more work. The strategy (essentially taken from \cite{km}), is the following.\\

In any locally convex space, there is a natural notion of bounded set. The collection of these bounded sets (``the von Neumann bornology") is not very sensitive to the locally convex topology: a classical theorem of Mackey in functional analysis shows that if one varies the topology while keeping the same dual space, the bounded sets remain the same. As a consequence, one can view the bounded sets from the perspective of the weak topology instead of the given topology. In the weak topology, it is natural and immediate that a subset of $E$ whose image by every linear functional is bounded must be itself bounded.\\  

On the other hand, for curves, being $\mathcal{C}^{\infty}$ is ultimately a bornological concept: the $\mathcal{C}^{\infty}$ curves remain the same if one changes the locally convex topology, while keeping the same underlying bornology. This follows from the fact that a $\mathcal{C}^1$ curve is locally Lipschitz (by the mean value theorem), and the Lipschitz condition (which is essentially bornological) implies continuity. Translating smoothness in terms of Lipschitz conditions (involving bounded sets), it becomes possible to use the dual characterization of boundedness given by Mackey's theorem, to obtain a dual characterization of smoothness.\\

In what follows, we recall, for the convenience of the reader and for later reference, the details of the above arguments, starting with Mackey's theorem, the cornerstone in proving Theorem \ref{weak smoothness of  curves} as well as other results in this paper. For a proof of the former theorem, see, e.g., Theorem 36.2 in \cite{tre} or Theorem 8.3.4
in  \cite{jar}. \\

\begin{thm} (Mackey's theorem)
Let $E$ be a locally convex space, and $B$ a subset of $E$. If $\ell(B)$ is bounded for every $\ell\in E'$, then $B$ is bounded.\\
\end{thm}

\begin{defi} Let $E$ be a locally convex space, and $c:\mathbb{R}\longrightarrow E$ a curve. 
\begin{enumerate} 
\item If $J$ is an open subset of $\mathbb{R}$, we say that $c$ is {\bf Lipschitz} on $J$ if the 
$ \displaystyle  \hbox{set}\,\,  \Big\{\frac{c(t_2)-c(t_1)}{t_2-t_1}\;;\;t_1,t_2\in J\hbox{ and }t_1\neq t_2\Big\}$ is bounded in $E$. 
\item We say that $c$ is {\bf locally Lipschitz} if every point in $\mathbb{R}$ has a neighborhood on which $c$ is Lipschitz. \\
\end{enumerate}
\end{defi}

\begin{defi}
Let $E$ be a locally convex space, and $c:\mathbb{R}\longrightarrow E$ a curve. For $k\in\mathbb{N}$, we say that $c$ is {\bf of class $\mathcal{L}ip^k$} if all the derivatives of $c$ up to order $k$ exist, and $c^{(k)}:\mathbb{R}\longrightarrow E$ is locally Lipschitz.\\
\end{defi}

If $A$ is any subset of $E$, we will denote by $\langle A\rangle$ the absolute convex hull of the closure of $A$. We will need the following version of the mean value theorem, for curves in a locally convex space.\\

\begin{thm} (Mean Value Theorem)
Let $E$ be a locally convex space, and $c:[a,b]\longrightarrow E$ a curve which is continuous on $[a,b]$, and differentiable on $]a,b[$. Then $$\frac{c(b)-c(a)}{b-a}\in\;\langle\{c'(t)\;;\;a<t<b\}\rangle$$
\end{thm} 

\noindent {\bf Proof.~} {\it cf.} \cite{km}, I.1.4. $\hfill\square$\\

\begin{cor} 
Let $E$ be a locally convex space, and $c:\mathbb{R}\longrightarrow E$ a curve which is differentiable on an open interval $J\subset\mathbb{R}$. If $c'$ is bounded on $J$, then $c$ is Lipschitz on $J$. 
\end{cor}

\noindent {\bf Proof.~} 
Let $t_1,t_2\in J$ with $t_1\neq t_2$. By the mean value theorem, $$\frac{c(t_2)-c(t_1)}{t_2-t_1}\in\;\langle\{c'(t)\;;\;t\in J\}\rangle$$
So $\Big\{\frac{c(t_2)-c(t_1)}{t_2-t_1}\;;\;t_1,t_2\in J\hbox{ and }t_1\neq t_2\Big\}\subset\langle\{c'(t)\;;\;t\in J\}\rangle$, and this last set is bounded since the absolute convex hull of any bounded set is bounded. $\hfill\square$\\

\begin{prop} \label{C^k and Lip^k}
Let $E$ be a locally convex space, and $c:\mathbb{R}\longrightarrow E$ a curve.
\begin{enumerate}
\item If $c$ is of class $\mathcal{C}^k$, then $c$ is of class $\mathcal{L}ip^{k-1}$.
\item If $c$ is of class $\mathcal{L}ip^{k-1}$, then $c$ is of class $\mathcal{C}^{k-1}$.
\end{enumerate}
\end{prop}

\noindent {\bf Proof.~} 1. If $c$ is $\mathcal{C}^k$, then $c^{(k-1)}$ is $\mathcal{C}^1$. This implies (via the preceding corollary) that $c^{(k-1)}$ is locally Lipschitz. Thus, $c$ is $\mathcal{L}ip^{k-1}$.\\
2. If $c$ is $\mathcal{L}ip^{k-1}$, then $c^{(k-1)}$ is locally Lipschitz. This implies that $c^{(k-1)}$ is continuous. Thus, $c$ is $\mathcal{C}^{k-1}$. $\hfill\square$\\

\begin{cor} 
Let $E$ be a locally convex space, and $c:\mathbb{R}\longrightarrow E$ a curve. Then $c$ is of class $\mathcal{C}^{\infty}$ if and only if $c$ is of class $\mathcal{L}ip^{k}$ for all $k\in\mathbb{N}$.
\end{cor}

\noindent {\bf Proof.~} If $c$ is smooth, then $c$ is $\mathcal{C}^k$ for every $k$. The first part of Proposition \ref{C^k and Lip^k} implies then that $c$ is $\mathcal{L}ip^{k-1}$ for every $k$. Conversely, if $c$ is $\mathcal{L}ip^{k}$ for every $k$, then the second part of Proposition \ref{C^k and Lip^k} implies that $c$ is $\mathcal{C}^{k}$ for every $k$, and so $c$ is smooth. $\hfill\square$\\

\begin{defi}
Let $E$ be a locally convex space. A curve $c:\mathbb{R}\longrightarrow E$ is said to be {\bf weakly smooth} if $\ell\circ c$ is smooth for every $\ell\in E'$.\\
\end{defi}

\begin{defi} Let $E$ be a locally convex space.
\begin{enumerate}
\item A sequence $(x_n)$ in $E$ is said to be {\bf Mackey-convergent} to a point $x\in E$ if there exists an absolutely convex bounded set $B\subset E$, and a sequence $(\mu_n)$ of real numbers converging to 0, such that $x_n-x\in \mu_n B$ for all $n$. 
\item A sequence $(x_n)$ in $E$ is said to be {\bf Mackey-Cauchy} if there exists an absolutely convex bounded set $B\subset E$, and a double sequence $(\mu_{n,m})$ of real numbers converging to 0, such that $x_n-x_m\in \mu_{n,m} B$ for all $n,m$. \\

We also have the same notions for nets (just replace ``sequence" by ``net" everywhere in the preceding definition).\\

\item $E$ is called {\bf Mackey-complete} if every Mackey-Cauchy net in $E$ converges.\\
\end{enumerate}
\end{defi}

\begin{rem} In fact, $E$ is Mackey-complete if and only if every Mackey-Cauchy sequence in $E$ converges, cf. \cite{km}, I.2.2. Note that every complete locally convex space is sequentially complete, and in turn sequential completeness implies Mackey-completeness.
For metrizable locally convex spaces the three notions of completeness coincide. Given a locally convex space $E$, the completion $\hat{E}$
yields a complete (and thus Mackey-complete) locally convex space together with a continuous linear embedding  $j:E\hookrightarrow \hat{E}$ 
having dense image. (Compare, e.g. Theorem 5.2 in \cite{tre}.) A continuous curve $c: \mathbb{R} \longrightarrow E$ yields then a continuous curve 
$\hat{c}= j\circ c:\mathbb{R}\longrightarrow \hat{E}$, and $c$ is weakly smooth if and only if $\hat{c}$ is weakly smooth.\\
\end{rem}

\begin{thm} \label{weak smoothness of  curves}
Let $E$ be a Mackey-complete locally convex space and $c:\mathbb{R}\longrightarrow E$ a curve. Then $c$ is smooth if and only if it is weakly smooth.
\end{thm}

\noindent {\bf Proof.~} If $c$ is smooth, then it is clear by chain rule that $\ell\circ c$ is smooth for every $\ell\in E'$. To prove the converse, suppose that for every $\ell\in E'$, the function $\ell\circ c$ is smooth. We first prove that $c$ is differentiable. Let $t_0\in\mathbb{R}$, and $J$ a compact interval about $t_0$. Set $\displaystyle q(t):=\frac{c(t)-c(t_0)}{t-t_0}$ for all $t\in J-\{t_0\}$. We need to show that $q$ has a limit as $t\to t_0$. Let $B:=\Big\{\frac{q(t)-q(t')}{t-t'}\;;\;t,t'\in J-\{t_0\}\hbox{ and }t\neq t'\Big\}$. Then 
$\ell(B)=\Big\{\frac{(\ell\circ q)(t)-(\ell\circ q)(t')}{t-t'}\;;\;t,t'\in J-\{t_0\}\hbox{ and }t\neq t'\Big\}$. Now for each $t\in J-\{t_0\}$, we have\\
 $\displaystyle (\ell\circ q)(t)=\frac{(\ell\circ c)(t)-(\ell\circ c)(t_0)}{t-t_0}$. Since $\ell\circ c$ is smooth, the same must be true for $\ell\circ q$. In particular, the mean value theorem implies the existence of $\tau\in \;]0,1[$ such that $(\ell\circ q)(t)-(\ell\circ q)(t')=(t-t')\;(\ell\circ q)'(t'+\tau(t-t'))$. Now $t'+s(t-t')\in J$ for every $s\in\;]0,1[$ (by convexity of $J$). Since $(\ell\circ q)'$ is continuous and $J$ is compact, we deduce that there is $M>0$ such that $|(\ell\circ q)'(t'+s(t-t'))|\leq M$ for all $s \in\; ]0,1[$. But then, $|(\ell\circ q)(t)-(\ell\circ q)(t')|\leq M|t-t'|$. Thus, $\ell(B)$ is bounded. Since $\ell$ was arbitrary, we deduce by Mackey's theorem that $B$ is bounded as well, i.e. $q$ is Lipschitz on $J-\{t_0\}$. Otherwise stated, the net $(q_t)_{t\in J-\{t_0\}}$ with $q_t:=q(t)$ is Mackey-Cauchy. Mackey-completeness of $E$ implies that $q$ has a continuous extension to $J$, which allows us to define $c'(t_0)$ as $q(t_0)$. Thus, $c$ is differentiable, and for every $\ell\in E'$, we have $(\ell\circ c)'=\ell\circ c'$ (by chain rule and linearity of $\ell$).\\ 
 
  Second, we prove that $c'$ is locally Lipschitz.  Let $t_0\in\mathbb{R}$, and $J$ a compact interval about $t_0$. Then for every $\ell\in E'$, the function $(\ell\circ c)'$ is Lipschitz on $J$ (since $(\ell\circ c)''$, being continuous, is bounded on $J$). Let $\displaystyle B_1=\Big\{\frac{c'(t_2)-c'(t_1)}{t_2-t_1}\;;\;t_1,t_2\in J\hbox{ and }t_1\neq t_2\Big\}$. Then, since $\ell\circ c'=(\ell\circ c)'$, we have $\displaystyle \ell(B_1)=\Big\{\frac{(\ell\circ c)'(t_2)-(\ell\circ c)'(t_1)}{t_2-t_1}\;;\;t_1,t_2\in J\hbox{ and }t_1\neq t_2\Big\}$. Since $(\ell\circ c)'$ is Lipschitz on $J$, we have that $\ell(B_1)$ is bounded. But $\ell$ is arbitrary. Using again Mackey's theorem, we deduce that $B_1$ is bounded, hence $c'$ is Lipschitz on $J$, and so $c'$ is locally Lipschitz, that is, $c$ is $\mathcal{L}ip^1$.\\
  
Replacing $c$ by $c'$ in the above chain of arguments, we arrive at the conclusion that $c'$ is  $\mathcal{L}ip^1$ as well, i.e. $c$ is $\mathcal{L}ip^2$. By induction, we conclude that $c$ is $\mathcal{L}ip^k$ for all $k$, i.e. $c$ is smooth. $\hfill\square$\\

\section{Smoothness of maps}

In this section, we start by recalling the notion of a smooth map from an open subset $U$ of a locally convex space $E$ into a locally convex space $F$. Among the various existing notions of (differentiable and) smooth maps, we choose what is called $\mathcal{C}^{\infty}_c$ in the book of Keller (\cite{kel}), going back at least to Bastiani (compare Definitions II.3.1 and II.2.2 in \cite{bas}). This definition of smoothness 
turns out to be appropriate for the construction of an applicable theory of infinite-dimensional manifolds and Lie groups, compare notably 
\cite{ham}, \cite{mil}, \cite{hw}, \cite{nee}.
The goal of this section is to generalize Theorem \ref{weak smoothness of  curves} to maps having infinite-dimensional source. Namely, we want to show that if $E$ is for example a Fr\'echet space or the continuous dual of a Fr\'echet-Schwartz space ({\it cf.} below for details), a map $\Phi:E \supset U\longrightarrow F$ is smooth if and only if it is scalarwise smooth. By {\it scalarwise smooth}, we mean that $f\circ\Phi\in\mathcal{C}^{\infty}(U,\mathbb{R})$ for every $f\in\mathcal{C}^{\infty}(F,\mathbb{R})$.\\

 As a matter of fact, this result is very easy to prove if we replace the notion of smoothness of Bastiani et al. by another one: the ``convenient smoothness" considered by Fr\"olicher, Kriegl and Michor (compare notably \cite{km}). A map $\Phi:U\longrightarrow F$ is said to be {\it conveniently smooth} if it sends smooth curves to smooth curves. Since the notion of convenient smoothness relies on curves, Theorem \ref{weak smoothness of  curves} immediately yields the equivalence between convenient smoothness and scalarwise convenient smoothness of maps.\\
 
For applications, it is thus highly desirable to establish that convenient smoothness coincides with smoothness in the sense of Bastiani-Hamilton-Milnor-Neeb, for relevant classes of locally convex spaces (``generalized Boman theorem"). That this is indeed the case was already observed in \cite{nee} in case where the source $E$ is a Fr\'echet space. Below we will first prove a more general result of this type, before giving our crucial characterization of smoothness in terms of weak resp. scalarwise smoothness.\\ 
 
Let $E$ and $F$ be locally convex spaces, and $U$ an open subset of $E$. For a continuous map $\Phi:U\longrightarrow F$, the {\it G\^{a}teaux derivative} of $\Phi$ at a point $x\in U$ in the direction of a vector $v\in E$ is defined by $$d\Phi_{|x}(v)=\lim_{t\to 0}\frac{\Phi(x+tv)-\Phi(x)}{t}$$ provided the limit exists.\\

\begin{defi} \label{BHMN}
The map $\Phi$ is said to be {\bf of class $\mathcal{C}^1$} if $d\Phi:U\times E\longrightarrow F$, $(x,v)\mapsto (d\Phi)(x,v):=d\Phi_{|x}(v)$ exists and is continuous. We define inductively a map $\Phi$ to be {\bf of class $\mathcal{C}^{k+1}$} if it is of class of class $\mathcal{C}^1$ and $d\Phi$ is
of class $\mathcal{C}^{k}$. Furthermore, a map is called being {\bf of class $\mathcal{C}^{\infty}$} or {\bf smooth} if it is $\mathcal{C}^{k}$ for all $k\in\mathbb N$.\\
\end{defi}

\begin{rem} If $E$ and $F$ are Banach spaces, being $\mathcal{C}^1$ in the above sense is weaker than the usual notion of $\mathcal{C}^1$ in the sense of Fr\'echet differentiability, which requires the map $x\mapsto d\Phi_{|x}$ to be continuous as map from $U$ to $\mathcal{L}(E,F)$, equipped with the operator norm topology. However, $\mathcal{C}^2$ in the above sense implies $\mathcal{C}^1$ in the usual Fr\'echet differentiability sense ({\it cf.} Proposition 2.7.1 in \cite{kel}, or Theorem I.7 in \cite{hw}), so that in Banach spaces, the two definitions lead to the same smooth maps.\\
\end{rem}

Next we turn to the notion of convenient smoothness. Before we recall its definition and main properties, let us already note that in general, it is possible to find conveniently smooth maps which are not even continuous! This hints to the fact that there should be a different topology which is more adapted to convenient smoothness, and for which conveniently smooth maps are automatically continuous. Since the definition of convenient smoothness relies on smooth curves, it is natural to use smooth curves to define this topology.\\

\begin{defi}
Let $E$ be a locally convex space. The {\bf $\hbox{\upshape{c}}^{\infty}$-topology} (also called {\bf Mackey-closure topology}) is the finest topology on $E$ making all the smooth curves $c:\mathbb{R}\longrightarrow E$ continuous. Subsets of $E$ open in this topology will be called {\bf $\hbox{\upshape{c}}^{\infty}$-open}.\\
\end{defi}

\begin{rem} The $\hbox{\upshape{c}}^{\infty}$-topology on $E$ is clearly finer than the given locally convex topology. Note that if $E$ is a Fr\'echet space or the dual of a Fr\'echet-Schwartz space, then the two topologies coincide (\cite{km}, Theorem I.4.11). \\ 

Note also that there are locally convex spaces such that the $\hbox{\upshape{c}}^{\infty}$-topology is not even a vector space topology and thus {\it a fortiori} does not equal the initially given topology. Examples of this phenomenon are strict inductive limits of strictly increasing sequences of infinite-dimensional Fr\'echet spaces, as e.g. $\mathcal{D}(M)$, the space of compactly supported smooth functions on a noncompact finite-dimensional smooth manifold. (Compare \cite{km}, Proposition I.4.26.)\\
\end{rem}

Let us now recall the precise definition and some of the properties of conveniently smooth maps. ({\it cf.} notably \cite{km}).\\

\begin{defi}
Let $E$ and $F$ be locally convex spaces, $U$ a $\hbox{\upshape{c}}^{\infty}$-open subset of $E$, and $\Phi:U\longrightarrow F$ a map. We say that $\Phi$ is {\bf conveniently smooth} if for every smooth curve $c:\mathbb{R}\longrightarrow U$, the curve $\Phi\circ c:\mathbb{R}\longrightarrow F$ is smooth. \\
\end{defi}

\begin{rem} It is easy to see that the composition of two conveniently smooth maps is conveniently smooth.\\
\end{rem}

\begin{prop} \label{convenient smoothness implies c-infinity continuity}

Let $E$ and $F$ be locally convex spaces, and $U$ an open subset of $E$. Any conveniently smooth map $\Phi:U\longrightarrow F$ is continuous when $E$ is equipped with the $\hbox{\upshape{c}}^{\infty}$-topology.
\end{prop}

\noindent {\bf Proof.~} Let $V$ be an open subset of $F$. To show that $\Phi^{-1}(V)$ is $\hbox{\upshape{c}}^{\infty}$-open in $U$, we need to show that $c^{-1}(\Phi^{-1}(V))$ is open in $\mathbb{R}$ for every $c\in\mathcal{C}^{\infty}(\mathbb{R},U)$. But $c^{-1}(\Phi^{-1}(V))=(\Phi\circ c)^{-1}(V)$, which is clearly open in $\mathbb{R}$ since the curve $\Phi\circ c$ is smooth, and therefore continuous. $\hfill\square$\\

We denote by $\mathcal{C}^{\infty}_{conv}(U,F)$ the set of conveniently smooth maps from $U$ to $F$. \\

The equivalence between smoothness and convenient smoothness, when it holds, is not trivial to establish, even for functions from $\mathbb{R}^d$ to $\mathbb{R}$. In this case, it was first proved by Boman in 1967 \cite{bom}.\\

\begin{thm} (Boman's theorem)
$$\mathcal{C}^{\infty}_{conv}(\mathbb{R}^2,\mathbb{R})=\mathcal{C}^{\infty}(\mathbb{R}^2,\mathbb{R}) \, .$$
\end{thm}

\noindent {\bf Proof.~} {\it cf.} \cite{bom} or \cite{km}, I.3.4. $\hfill\square$\\

Now we concentrate on the structure of $\mathcal{C}^{\infty}_{conv}(U,F)$.\\

\begin{prop} 
Let $E$ and $F$ be locally convex spaces, $U$ a $\hbox{\upshape{c}}^{\infty}$-open subset of $E$. Then $\mathcal{C}^{\infty}_{conv}(U,F)$ is a locally convex space for the coarsest topology making the maps $c^*:\mathcal{C}^{\infty}_{conv}(U,F)\longrightarrow \mathcal{C}^{\infty}(\mathbb{R},F)$, $\Phi\mapsto c^*\Phi=\Phi\circ c$ for all $c\in \mathcal{C}^{\infty}(\mathbb{R},U)$ continuous.
\end{prop}

\noindent {\bf Proof.~} It is not difficult to check that this topology on $\mathcal{C}^{\infty}_{conv}(U,F)$ is defined by the family of seminorms $\mathcal{P}=\{p_{K,\alpha,q,c}\;;\;K\hbox{ compact in $\mathbb{R}$},\\ \alpha\in\mathbb{N}, q\hbox{ cont. seminorm on $F$}, c\in\mathcal{C}^{\infty}(\mathbb{R},U)\}$, where we set $\displaystyle p_{K,\alpha,q,c}(\Phi):=\sup_{t\in K}q((\Phi\circ c)^{(\alpha)}(t))$
for $\Phi\in\mathcal{C}^{\infty}_{conv}(U,F)$. Moreover, if $\Phi\in\mathcal{C}^{\infty}_{conv}(U,F)-\{0\}$, let $x\in U$ be such that $\Phi(x)\neq 0$, and $c:=k_x:\mathbb{R}\longrightarrow U$ the constant curve at $x$. There exists a continuous seminorm $q$ on $F$ such that $q(\Phi(x))\neq 0$. Take $\alpha=0$ and $K=\{0\}$. Then $\displaystyle p_{K,\alpha,q,c}(\Phi)=\sup_{t\in\{0\}}q((\Phi\circ c)(t))=q((\Phi\circ c)(0))=q(\Phi(x))\neq 0$. Thus, $\mathcal{P}$ is separating. $\hfill\square$\\

One of the major benefits of working with conveniently smooth maps is the fact that they give rise to a Cartesian closed category, as the next theorem will show.\\

\begin{thm} \label{Cartesian closedness}
Let $E_1$, $E_2$ and $F$ be locally convex spaces, $U_1$ and $U_2$ $\hbox{\upshape{c}}^{\infty}$-open subsets of $E_1$ and $E_2$ respectively. Then, as sets,
$$\mathcal{C}^{\infty}_{conv}(U_1\times U_2,F)\cong\mathcal{C}^{\infty}_{conv}(U_1,\mathcal{C}^{\infty}_{conv}(U_2,F))\,.$$
\end{thm}

\noindent {\bf Proof.~} {\it cf.}   Theorem I.3.12 in \cite{km}. $\hfill\square$\\

A first consequence of Cartesian closedness is the following generalization of Boman's theorem.\\

\begin{cor} \label{conveniently smooth on R^n = smooth on R^n}
Let $E$ be a locally convex space. Then,
$$\mathcal{C}^{\infty}_{conv}(\mathbb{R}^n,E)=\mathcal{C}^{\infty}(\mathbb{R}^n,E)$$
\end{cor}

\noindent {\bf Proof.~} Let $\Phi\in\mathcal{C}^{\infty}_{conv}(\mathbb{R}^n,E)$. By Theorem \ref{Cartesian closedness}, we have 
$\mathcal{C}^{\infty}_{conv}(\mathbb{R}^n,E)\cong$ $ \mathcal{C}^{\infty}_{conv}(\mathbb{R}^{n-1},\mathcal{C}^{\infty}(\mathbb{R},E))$. In particular, for every $x_1,...,x_{i-1},x_{i+1},...,x_n\in\mathbb{R}$, the partial map $\Phi(x_1,...,\bullet,...,x_n)$, $y\mapsto \Phi(x_1,...,x_{i-1},y,x_{i+1},...,x_n)$ lies in $\mathcal{C}^{\infty}(\mathbb{R},E)$. Then $\displaystyle\frac{\partial\Phi}{\partial x_i}(x_1,...,x_i,...,x_n)=$\\
$\displaystyle \lim_{t \to 0}\frac{\Phi(x_1,...,x_i+t,...,x_n)-\Phi(x_1,...,x_i,...,x_n)}{t}=$\\
$\displaystyle {\lim_{t \to 0}}\frac{\Phi(x_1,...,\bullet,...,x_n)(x_i+t)-\Phi(x_1,...,\bullet,...,x_n)(x_i)}{t}{=}\Phi(x_1,...,\bullet,...,x_n)'(x_i)$. Thus, all first-order partial derivatives of $\Phi$ exist. Inductively, one obtains the existence of all higher-order partial derivatives of $\Phi$. This proves that $\Phi\in \mathcal{C}^{\infty}(\mathbb{R}^n,E)$. Conversely, if $\Phi\in\mathcal{C}^{\infty}(\mathbb{R}^n,E)$, then $\Phi\in\mathcal{C}^{\infty}_{conv}(\mathbb{R}^n,E)$ by the chain rule. $\hfill\square$\\

Another consequence of Cartesian closedness is the convenient smoothness of the differential.\\

\begin{prop}
Let $E$ and $F$ be locally convex spaces, $U$ a $\hbox{\upshape{c}}^{\infty}$-open subset of $E$, and $\Phi\in \mathcal{C}^{\infty}_{conv}(U,F)$. For every $x\in U$ and $v\in E$, set $\displaystyle d\Phi_{|x}(v):=\lim_{t\to 0}\frac{\Phi(x+tv)-\Phi(x)}{t}$. Then $d\Phi\in\mathcal{C}^{\infty}_{conv}(U\times E,F)$. 
\end{prop}

\noindent {\bf Proof.~} We claim that the map $\delta:\mathcal{C}^{\infty}_{conv}(U,F)\times U\times E\longrightarrow F$ defined by $$\delta(\Phi,x,v):=d\Phi_{|x}(v)=\lim_{s\to 0}\frac{\Phi(x+sv)-\Phi(x)}{s}$$ is conveniently smooth. Indeed, let $c=(\tilde{\Phi},\tilde{x},\tilde{v}):\mathbb{R}\longrightarrow \mathcal{C}^{\infty}_{conv}(U,F)\times U\times E$ be a smooth curve, and set $h(t,s):=\tilde{\Phi}(t)(\tilde{x}(t)+s\tilde{v}(t))$. Then $\displaystyle(\delta\circ c)(t)=\delta(c(t))=\delta(\tilde{\Phi}(t),\tilde{x}(t),\tilde{v}(t))=$\\
$\displaystyle \lim_{s\to 0}\frac{\tilde{\Phi}(t)(\tilde{x}(t)+s\tilde{v}(t))-\tilde{\Phi}(t)(\tilde{x}(t))}{s}=\lim_{s\to 0}\frac{h(t,s)-h(t,0)}{s}=\frac{\partial h}{\partial s}(t,0)$. Since $h$ is clearly conveniently smooth (and smooth by Corollary \ref{conveniently smooth on R^n = smooth on R^n}), we deduce that the curve $\delta\circ c$ is smooth. Thus, $\delta\in\mathcal{C}^{\infty}_{conv}(\mathcal{C}^{\infty}_{conv}(U,F)\times U\times E\;,\;F)$, and therefore $d:=\hat{\delta}\in \mathcal{C}^{\infty}_{conv}(\mathcal{C}^{\infty}_{conv}(U,F)\;,\;\mathcal{C}^{\infty}_{conv}(U\times E,F))$. In particular, $d\Phi\in\mathcal{C}^{\infty}_{conv}(U\times E,F)$. $\hfill\square$\\

An immediate consequence of the above proposition, already observed in \cite{nee} for Fr\'echet spaces, is a further generalization of Boman's theorem, stating that for a reasonable class of locally convex spaces, convenient smoothness coincides with smoothness (in the sense of 
the above definition).\\

\begin{prop}
Let $E$ and $F$ be locally convex spaces, and $U$ an open subset of $E$. Assume that for every $n\geq 1$, the $\hbox{\upshape{c}}^{\infty}$-topology on $E^n$ is the same as the product topology for the given topology on $E$. Then $$\mathcal{C}^{\infty}_{conv}(U,F)=\mathcal{C}^{\infty}(U,F)$$
\end{prop}

\noindent {\bf Proof.~} By the chain rule, smoothness implies convenient smoothness. For the nontrivial direction, suppose $\Phi:U\longrightarrow F$ is conveniently smooth. By the preceding proposition, $d\Phi:U\times E\longrightarrow F$ is conveniently smooth as well. Proposition \ref{convenient smoothness implies c-infinity continuity} implies then that $d\Phi$ is continuous for the $\hbox{\upshape{c}}^{\infty}$-topology, therefore continuous by our assumption on $E$. Thus, $\Phi$ is $\mathcal{C}^1$, and by induction, one obtains that $\Phi$ is smooth. $\hfill\square$\\

We are ready to state and prove the main result of this section.\\

\begin{thm} \label{scalarwise smoothness implies smoothness}
Let $E$ and $F$ be locally convex space spaces, and $U$ an open subset of $E$. Assume that for every $n\geq 1$, the $\hbox{\upshape{c}}^{\infty}$-topology on $E^n$ is the same as the product topology for the given topology on $E$ and that F is Mackey-complete. Let $\Phi:U\longrightarrow F$ be a continuous map. Then the following are equivalent:
\begin{enumerate}
\item [(i)] $\Phi$ is smooth.
\item [(ii)] For every open subset $V$ of $F$ containing $\Phi(U)$ and every $f\in\mathcal{C}^{\infty}(V,\mathbb{R})$, we have $f\circ\Phi\in \mathcal{C}^{\infty}(U,\mathbb{R})$.

\item [(iii)] For every $f\in \mathcal{C}^{\infty}(F,\mathbb{R})$, we have $f\circ\Phi\in \mathcal{C}^{\infty}(U,\mathbb{R})$ (i.e. $\Phi$ is scalarwise smooth).
\item [(iv)] For every $\ell\in F'$, we have $\ell\circ\Phi\in \mathcal{C}^{\infty}(U,\mathbb{R})$ (i.e. $\Phi$ is weakly smooth). 
\end{enumerate}
\end{thm}

\noindent {\bf Proof.~}
$(i)\Longrightarrow(ii)$ is evident by chain rule, and $(ii)\Longrightarrow(iii)$ follows upon taking $V=F$. $(iii)\Longrightarrow(iv)$ follows from the fact that every continuous linear map is smooth. It remains to show $(iv)\Longrightarrow(i)$. Suppose that for every $\ell\in F'$, we have $\ell\circ \Phi\in \mathcal{C}^{\infty}(U,\mathbb{R})$. For every smooth curve $c:\mathbb{R}\longrightarrow U$, the function $\ell\circ(\Phi\circ c)$ is smooth since $\ell\circ(\Phi\circ c)=(\ell\circ\Phi)\circ c$. By Theorem \ref{weak smoothness of curves}, the curve $\Phi\circ c:\mathbb{R}\longrightarrow F$ is smooth. This means that the map $\Phi:U\longrightarrow F$ is conveniently smooth. By the above proposition, we conclude that $\Phi$ must be smooth. $\hfill\square$\\

\section{Infinite-dimensional manifolds}

\begin{defi} A {\bf class of l.c. model spaces or l.c. models} is a subclass $\mathcal{E}$ of the class of real (Hausdorff) locally convex spaces such that $\mathcal{E}$ contains the numerical spaces $\mathbb{R}^n$ for all $n\in\mathbb{N}$.\\
\end{defi}
 
Typical examples are the class of numerical spaces $\mathbb{R}^n$ for all $n\in\mathbb{N}$, the class of finite-dimensional vector spaces, the class of Banach spaces, the class of Fr\'echet spaces, and the class of all locally convex spaces.\\
 
\begin{defi} \label{atlas} Let $\mathcal{E}$ be a class of l.c. model spaces, and $M_0$ a Hausdorff topological space. 
\begin{enumerate}
\item A {\bf smooth $\mathcal{E}$-atlas} on $M_0$ is a family of pairs $\mathcal{A}=\{(U_{\alpha},\varphi_{\alpha})\;;\;\alpha\in A\}$ such that $\{U_{\alpha}\;;\;\alpha\in A\}$ is an open cover of $M_0$ and for every $\alpha\in A$, there exists a space $E_{\alpha}$ from the class $\mathcal{E}$ and a homeomorphism $\varphi_{\alpha}:U_{\alpha}\longrightarrow\varphi_{\alpha}(U_{\alpha})\subset E_{\alpha}$ such that the following compatibility condition is satisfied: for every $\alpha,\beta\in A$ such that $U_{\alpha\beta}=U_{\alpha}\cap U_{\beta}\neq\emptyset$, the transition map $\varphi_{\alpha\beta}=\varphi_{\alpha}\circ\varphi_{\beta}^{-1}
:\varphi_{\beta}(U_{\alpha}\cap U_{\beta})\longrightarrow \varphi_{\alpha}(U_{\alpha}\cap U_{\beta})$ is smooth (in the sense of Definition \ref{BHMN}).
\item Given a smooth $\mathcal{E}$-atlas $\mathcal{A}$ on $M_0$,  a {\bf (compatible) chart} on $M_0$ is a pair $(U,\varphi)$ where $U$ is an open subset of $M_0$ and $\varphi$ is a homeomorphism from $U$ onto an open subset $\varphi(U)$ of some space from the class $\mathcal{E}$ such that
the transition map between $(U,\varphi)$ and every $(U_{\alpha},\varphi_{\alpha})\in\mathcal{A}$ is smooth.
\item A smooth $\mathcal{E}$-atlas $\mathcal{A}$ is said to be {\bf maximal} if any chart $(U,\varphi)$ compatible with $\mathcal{A}$ belongs already 
to $\mathcal{A}$. A maximal atlas is also called a {\bf smooth structure}, the class $\mathcal{E}$ being understood.\\
\end{enumerate}
\end{defi}

\begin{defi}
Let $\mathcal{E}$ be a class of l.c. model spaces. A {\bf smooth $\mathcal{E}$-manifold} is a pair $M=(M_0,\mathcal{A})$ where $M_0$ is a Hausdorff topological space and $\mathcal{A}$ is a maximal smooth $\mathcal{E}$-atlas on $M_0$.\\
\end{defi}  

\begin{rem}\
\begin{enumerate}
\item Especially in finite dimensions, $M_0$ is often required to be second countable. However, as we are mainly interested in the generalization to infinite dimensions, we do not insist on this condition here. Let us remark that we cannot insist on first countability either, since a non-metrizable locally convex space is not first countable (see e.g. \cite{mv}, Proposition 25.1).  
\item  Given a smooth $\mathcal{E}$-atlas $\mathcal{A}$, there is always a unique maximal smooth $\mathcal{E}$-atlas containing $\mathcal{A}$, obtained by adjoining to $\mathcal{A}$ all the charts that are compatible with $\mathcal{A}$. 
\item If $\mathcal{E}$ and $\mathcal{E}'$ are classes of l.c. models such that $\mathcal{E}\subset\mathcal{E}'$, then a  smooth $\mathcal{E}$-manifold is obviously a smooth $\mathcal{E}'$-manifold. Denoting the class of all locally convex spaces by $\mathbf{LCS}$, every  smooth $\mathcal{E}$-manifold is then a smooth $\mathbf{LCS}$-manifold.\\
\end{enumerate}
\end{rem}

\begin{defi}
Given classes  $\mathcal{E}$ and $\mathcal{F}$ of l.c. models, let $M=(M_0,\mathcal{A})$, resp. $N=(N_0,\mathcal{B})$ be a smooth $\mathcal{E}$-
resp. $\mathcal{F}$-manifold, and $\Phi:M_0\longrightarrow N_0$ a continuous map. We say that $\Phi$ is {\bf smooth} if the following condition is fulfilled: for every $(U,\varphi)\in\mathcal{A}$ and every $(V,\psi)\in\mathcal{B}$ such that $\Phi(U)\subset V$, the following map (between open sets in locally convex spaces) is smooth:
$$\psi\circ\Phi_{|U} \circ \varphi^{-1}:\varphi(U)\longrightarrow\psi(V)\,.$$
\end{defi}

\begin{rem}\
\begin{enumerate}
\item If $\mathcal{A}'\subset\mathcal{A}$ resp. $\mathcal{B}'\subset\mathcal{B}$ is a smooth atlas (not necessarily maximal) of $M$ resp. $N$, the above condition is equivalent to the following: for every point $p\in M$, there exists a chart $(U,\varphi)$ containing $p$ and compatible with $\mathcal{A}'$ and a chart $(V,\psi)$ compatible with $\mathcal{B}'$ such that $\Phi(U)\subset V$ and such that the map $\psi\circ\Phi_{|U} \circ \varphi^{-1}:\varphi(U)\longrightarrow\psi(V)$ is smooth.
\item We notably obtain the notion of a smooth function $f$ on $M$, namely a function $f:M_0\longrightarrow\mathbb{R}$ such that for every chart $(U,\varphi)\in \mathcal{A}'$ (atlas contained in $\mathcal{A}$), $f\circ \varphi^{-1}:\varphi(U)\longrightarrow\mathbb{R}$ is a smooth function on the open set $\varphi(U)$ (contained in some space from the class $\mathcal{E}$).
\item Since an open subset $U\subset M_0$ inherits obviously a smooth structure upon restricting the smooth structure to $U$, we have a notion of smooth map/function defined on $U$.\\
\end{enumerate}
\end{rem} 

\begin{defi}
Given a class of l.c. model spaces $\mathcal{E}$, let $M=(M_0,\mathcal{A})$ be a smooth $\mathcal{E}$-manifold. The 
{\bf sheaf $\mathcal{C}^{\infty}_M$ of smooth functions} on $M$ is the subsheaf of the sheaf $\mathcal{C}^0_{M_0}$ of 
continuous (real-valued) functions on the topological space $M_0$, defined as the contravariant functor 
$$\mathcal{C}^{\infty}_M:\mathbf{Open}(M_0)\longrightarrow \mathbb{R}\mathbf{Alg}_{com}$$
 assigning to every open subset $U$ of $M_0$ the commutative $\mathbb{R}$-algebra $\mathcal{C}^{\infty}_{M}(U)$ of smooth functions on $U$. Here, $\mathbf{Open}(M_0)$ is the category of open subsets of $M_0$ with inclusions as morphisms, and $\mathbb{R}\mathbf{Alg}_{com}$ is the category of unital commutative associative $\mathbb{R}$-algebras with unital $\mathbb{R}$-algebra homomorphisms as morphisms.\\
\end{defi}

\begin{rem} \label{manifolds are ringed}
Obviously, $(M_0,\mathcal{C}^{\infty}_M)$ is a ringed space, with the stalks $(\mathcal{C}^{\infty}_M)_p$ (for $p\in M_0$) being local unital $\mathbb{R}$-algebras with maximal ideals $\mathfrak{m}_p=\{ f_p\in(\mathcal{C}^{\infty}_M)_p\;|\;f(p)=0\}$. In short, $(M_0,\mathcal{C}^{\infty}_M)$ is a locally ringed space. Furthermore, $(M_0,\mathcal{C}^{\infty}_M)$ is, as a locally ringed space, locally isomorphic to models $(D_0,(\mathcal{C}^{\infty}_E)_{|D_0})$, where $E$ is a space from the class $\mathcal{E}$, and $D_0$ is an open subset of $E$. The atlas $\{(D_0,j_{D_0})\}$, where $j_{D_0}:D_0\hookrightarrow E$ is the canonical inclusion, is obviously smooth and thus contained in a maximal atlas $\mathcal{A}_{D_0}$. Calling $D$ the smooth manifold given by $(D_0,\mathcal{A}_{D_0})$, we obviously have $(\mathcal{C}^{\infty}_E)_{|D_0}=\mathcal{C}^{\infty}_D$.\\
\end{rem}

Recall that if $(X_0,\mathcal{O}_X)$ and $(Y_0,\mathcal{O}_Y)$ are locally ringed spaces, a {\it morphism of locally ringed spaces} between them is a pair $\Phi=(\Phi_0,\Phi^{\sharp})$, where $\Phi_0:X_0\longrightarrow Y_0$ is a continuous map and $\Phi^{\sharp}:\Phi_0^{-1}\mathcal{O}_Y\longrightarrow\mathcal{O}_X$ is a morphism of sheaves of unital $\mathbb{R}$-algebras such that for every $x\in X_0$, the induced unital $\mathbb{R}$-algebra homomorphism $\Phi^{\sharp}_x:(\mathcal{O}_Y)_{\Phi_0(x)}\longrightarrow(\mathcal{O}_X)_x$ is {\it local}, i.e. the image by $\Phi^{\sharp}_x$ of the maximal ideal of $(\mathcal{O}_Y)_{\Phi_0(x)}$ is contained in the maximal ideal of $(\mathcal{O}_X)_x$. Let us also recall that $\Phi^{\sharp}$ can be equivalently viewed as a sheaf morphism $\mathcal{O}_Y\longrightarrow(\Phi_0)_*\mathcal{O}_X$. In the sequel, we apply both formulations without further comment.\\

A locally ringed space $(X_0,\mathcal{O}_X)$ is here said to be {\it reduced} if $\mathcal{O}_X$ is a subsheaf of the sheaf of continuous functions $\mathcal{C}^0_{X_0}$. If $\Phi=(\Phi_0,\Phi^{\sharp}):(X_0,\mathcal{O}_X)\longrightarrow(Y_0,\mathcal{O}_Y)$ is a morphism between {\it reduced} locally ringed spaces, then $\Phi^{\sharp}$ is given by $\Phi_0^{*}$, the pullback by $\Phi_0$ (and so the morphism $\Phi$ is completely determined by the underlying continuous map $\Phi_0$). Consequently, we often write $\Phi$ instead of $\Phi_0$ when the spaces are reduced. Of course, we never do so in the non-reduced case, since $\Phi^{\sharp}$ is then part of the data.\\  

\begin{defi}
Let $\mathcal{E}$ be a class of l.c. models. A {\bf structure sheaf-smooth $\mathcal{E}$-manifold} is a reduced locally ringed space whose underlying topological space
is Hausdorff, and which, as a locally ringed space, is locally isomorphic to models $(D_0,\mathcal{C}^{\infty}_D)$, where $D_0$ is an open subset of a space from the class $\mathcal{E}$. \\
\end{defi}

\begin{rem} By Remark \ref{manifolds are ringed}, every smooth $\mathcal{E}$-manifold is, in a natural way, a structure sheaf-smooth $\mathcal{E}$-manifold. The converse is nontrivial. Using the preceding section, we can nevertheless show the following results.\\
\end{rem}

\begin{nota}
The class of l.c. models made of all the Mackey-complete locally convex spaces will be denoted by  $\mathcal{M}$. \\
\end{nota}

\begin{thm} \label{smooth equals ringed morphism}
Let $\mathcal{E}$ be the class of locally convex spaces $E$ such that for every $n\geq 1$, the $\hbox{\upshape{c}}^{\infty}$-topology on $E^n$ is the same as the product topology for the given topology on $E$. If $M=(M_0,\mathcal{A})$ is a smooth $\mathcal{E}$-manifold, $N=(N_0,\mathcal{B})$ a smooth $\mathcal{M}$-manifold, and $\Phi:M_0\longrightarrow N_0$ a continuous map, then $\Phi$ is smooth if and only if $\Phi:(M_0,\mathcal{C}^{\infty}_M)\longrightarrow (N_0,\mathcal{C}^{\infty}_N)$ is a morphism of locally ringed spaces.
\end{thm}

\noindent {\bf Proof.~} By the chain rule, smooth maps are morphisms of locally ringed spaces. Assume now that $\Phi$ is a morphism of locally ringed spaces. Let $(U,\varphi)$ and $(V,\psi)$ be charts of $M$ (resp. $N$) with values in $E\in\mathcal{E}$ 
(resp. $F\in \mathcal{M}$) such that $\Phi(U)\subset V$. The continuous map $\tilde{\Phi}:=\psi\circ\Phi_{|U} \circ \varphi^{-1}:E\supset\varphi(U)\longrightarrow\psi(V)\subset F$ is then a morphism of locally ringed spaces as well (with respect to the natural structure sheaves). Thus, $\tilde{\Phi}$ satisfies condition (ii) of Theorem \ref{scalarwise smoothness implies smoothness}, which in turn implies that $\tilde{\Phi}$ is smooth. It follows that $\Phi$ itself is smooth. $\hfill\square$\\

The above theorem has the following very important consequence.\\

\begin{cor} \label{manifolds = structure sheaf-smooth manifolds}
Let $\mathcal{E}$ be the class of Mackey-complete locally convex spaces $E$ such that for every $n\geq 1$, the $\hbox{\upshape{c}}^{\infty}$-topology on $E^n$ is the same as the product topology for the given topology on $E$. For every smooth $\mathcal{E}$-manifold in the structure-sheaf sense $(M_0,\mathcal{O}_M)$, there is a canonical maximal atlas $\mathcal{A}$ on $M_0$ such that $(M_0,\mathcal{A})$ is a smooth $\mathcal{E}$-manifold fulfilling $\mathcal{C}^{\infty}_M=\mathcal{O}_M$. Furthermore, the maximal atlas $\mathcal{A}$ is uniquely determined by the condition $\mathcal{C}^{\infty}_M=\mathcal{O}_M$.
\end{cor}

\noindent {\bf Proof.~} 
Let $\mathcal{A}'=\{(U_{\alpha},\varphi_{\alpha})\;;\;\alpha\in A\}$ be a family of pairs such that $\{U_{\alpha}\;;\;\alpha\in A\}$ is a covering of $M_0$ and for all $\alpha\in A$, $\varphi_{\alpha}:(U_{\alpha},(\mathcal{O}_M)_{|U_{\alpha}})\longrightarrow  (\varphi_{\alpha}(U_{\alpha}),(\mathcal{C}^{\infty}_{E_{\alpha}})_{|\varphi_{\alpha}(U_{\alpha})})$ (with $E_{\alpha}\in\mathcal{E}$ and $\varphi_{\alpha}(U_{\alpha})\subset E_{\alpha}$) is an isomorphism of locally ringed spaces. Furthermore, the continuous transition map $\varphi_{\alpha\beta}:\varphi_{\beta}(U_{\alpha\beta})\longrightarrow \varphi_{\alpha}(U_{\alpha\beta})$ is an isomorphism of the locally ringed spaces $(\varphi_{\beta}(U_{\alpha\beta}),(\mathcal{C}^{\infty}_{E_{\beta}})_{|\varphi_{\beta}(U_{\alpha\beta})})$ and $(\varphi_{\alpha}(U_{\alpha\beta}),(\mathcal{C}^{\infty}_{E_{\alpha}})_{|\varphi_{\alpha}(U_{\alpha\beta})})$. By the preceding theorem, the map $\varphi_{\alpha\beta}$ is then already smooth. It follows that $\mathcal{A}'$ is a smooth atlas on $M_0$. We observe that the unique maximal atlas $\mathcal{A}$ containing $\mathcal{A}'$ is canonically associated to the given locally ringed space. Moreover, the corresponding sheaf of smooth functions $\mathcal{C}^{\infty}_M$ is equal to the structure sheaf $\mathcal{O}_M$. Now assume that $\tilde{M}=(M_0,\mathcal{B})$ is a smooth $\mathcal{E}$-manifold such that $\mathcal{C}^{\infty}_{\tilde{M}}=\mathcal{O}_M$. Then $(\hbox{id}_{M_0},\hbox{id}_{M_0}^*):(M_0,\mathcal{C}^{\infty}_{M})\longrightarrow (M_0,\mathcal{C}^{\infty}_{\tilde{M}})$ is an isomorphism of locally ringed spaces. Using again the preceding theorem, we obtain that $\hbox{id}_{M_0}$ is a smooth diffeomorphism between $(M_0,\mathcal{A})$ and $(M_0,\mathcal{B})$. This, of course, implies that $\mathcal{A}$ and $\mathcal{B}$ are smoothly compatible atlases, and by maximality of $\mathcal{A}$ and $\mathcal{B}$, that $\mathcal{A}=\mathcal{B}$. $\hfill\square$\\

\begin{rem} The preceding theorem and its corollary show that for important classes of l.c. model spaces, as e.g. the class of Fr\'echet spaces, we can encode smoothness completely in sheaf-theoretic language. This approach simplifies the verification and application of smoothness in infinite dimensions, and allows generalizations to the non-reduced case, as e.g. for ``Fr\'echet supermanifolds". In general, i.e., for arbitrary classes of l.c. model spaces, resp. for our (regular) local models replaced by more general local models, as e.g., complex-analytic sets in open subsets of complex locally convex spaces, this encoding might not be possible anymore. In order to circumvent this problem, Douady introduced in \cite{dou} the notions of the next definition.\\
\end{rem}

\begin{rem} Let $\mathcal{C}$ be a category, and denote, for a topological space $Z_0$, the category of sheaves (of sets) on $Z_0$ by $\mathbf{Sh}_{Z_0}$. If $\Phi_0:X_0\longrightarrow Y_0$ is a continuous map of topological spaces and $\mathcal{G}:\mathcal{C}\longrightarrow \mathbf{Sh}_{Y_0}$ is a covariant functor, then $\Phi_0^{-1}\mathcal{G}:\mathcal{C}\longrightarrow \mathbf{Sh}_{X_0}$ defined by $(\Phi_0^{-1}\mathcal{G})(A)=\Phi_0^{-1}(\mathcal{G}(A))$ is a covariant functor, called the {\bf inverse image of $\mathcal{G}$ by $\Phi_0$}.\\
\end{rem}

\begin{defi} Let $\mathcal{C}$ be a category.
\begin{enumerate}
\item A {\bf $\mathcal{C}$-functored space} is a pair $(X_0,\mathcal{O}_X^{\mathcal{C}})$ where $X_0$ is a topological space, and $\mathcal{O}_X^{\mathcal{C}}:\mathcal{C}\longrightarrow \mathbf{Sh}_{X_0}$ is a covariant functor from $\mathcal{C}$ to the category of sheaves on $X_0$. We say that $\mathcal{O}_X^{\mathcal{C}}$ is the {\bf structure functor} of $(X_0,\mathcal{O}_X^{\mathcal{C}})$.
\item If $(X_0,\mathcal{O}^{\mathcal{C}}_X)$ and $(Y_0,\mathcal{O}^{\mathcal{C}}_Y)$ are $\mathcal{C}$-functored spaces, a {\bf morphism of $\mathcal{C}$-functored spaces} between them is a pair $\Phi=(\Phi_0,\Phi^{\sharp})$ where $\Phi_0:X_0\longrightarrow Y_0$ is a continuous map and $\Phi^{\sharp}:(\Phi_0)^{-1}\mathcal{O}^{\mathcal{C}}_Y\longrightarrow\mathcal{O}^{\mathcal{C}}_X$ is a natural transformation between the two functors $\Phi_0^{-1}\mathcal{O}^{\mathcal{C}}_Y:\mathcal{C}\longrightarrow\mathbf{Sh}_{X_0}$ and $\mathcal{O}^{\mathcal{C}}_X:\mathcal{C}\longrightarrow\mathbf{Sh}_{X_0}$.\\
\end{enumerate}
\end{defi}

\begin{rem}
 Given a class $\mathcal{F}$ of l.c. models, we continue to use the symbol $\mathcal{F}$, by a slight abuse of notation, to denote the category whose objects are open subsets of spaces from the class $\mathcal{F}$, and whose morphisms are the smooth maps between them (in the sense of Definition \ref{BHMN}). Such a category $\mathcal{F}$ will be referred to as a ``category of l.c. models".\\
\end{rem}

\begin{defi}
Given a category of l.c. models $\mathcal{F}$, an $\mathcal{F}$-functored space $(X_0,\mathcal{O}_X^{\mathcal{F}})$ is said to be {\bf reduced} if for every $V\in\hbox{\upshape{Ob}}(\mathcal{F})$ and every open set $U\subset X_0$, we have $\mathcal{O}_X^{\mathcal{F}}(V)(U)\subset \mathcal{C}^0(U,V)$, and $\mathcal{O}_X:=\mathcal{O}_X^{\mathcal{F}}(\mathbb{R})$ is a sheaf having local unital $\mathbb{R}$-algebras as stalks.\\
\end{defi}

\begin{rem} If $(X_0,\mathcal{O}_X^{\mathcal{F}})$ is a reduced $\mathcal{F}$-functored space, then $(X_0,\mathcal{O}_X)$ is a reduced locally ringed space.\\
\end{rem}

\begin{defi}
Given a class of l.c. models $\mathcal{E}$, let $M=(M_0,\mathcal{A})$ be a smooth $\mathcal{E}$-manifold. For any category of l.c. models $\mathcal{F}$, consider the functor $(\mathcal{C}^{\infty})^{\mathcal{F}}_M:\mathcal{F}\longrightarrow \mathbf{Sh}_{M_0}$ defined as follows: for every $V\in\hbox{\upshape{Ob}}(\mathcal{F})$ and every open set $U\subset M_0$, let $(\mathcal{C}^{\infty})^{\mathcal{F}}_M(V)(U):=\mathcal{C}^{\infty}(U,V)$. Then $(M_0,(\mathcal{C}^{\infty})^{\mathcal{F}}_M)$ is called the {\bf canonical $\mathcal{F}$-functored space associated with $M$}. \\
\end{defi}

\begin{rem} \label{manifolds are functored}
Obviously, for any smooth $\mathcal{E}$-manifold $M$, $(M_0,(\mathcal{C}^{\infty})^{\mathcal{F}}_M)$ is a reduced $\mathcal{F}$-functored space. Furthermore, $(M_0,(\mathcal{C}^{\infty})^{\mathcal{F}}_M)$ is, as a $\mathcal{F}$-functored space, locally isomorphic to model $\mathcal{F}$-functored spaces $(D_0,(\mathcal{C}^{\infty})^{\mathcal{F}}_D)$, where $D_0$ is an open subset of a space from the class $\mathcal{E}$. \\
\end{rem}

\begin{defi}
Let $\mathcal{E}$ be a class of l.c. models. 
A {\bf structure functor-smooth $\mathcal{E}$-manifold} is a reduced $\bf{LCS}$-functored space 
$(M_0,\mathcal{O}_M^{\mathbf{LCS}})$ whose underlying topological space $M_0$ is Hausdorff, and
which, as a $\mathbf{LCS}$-functored space, is locally isomorphic to model functored spaces $(D_0,(\mathcal{C}^{\infty})^{\mathbf{LCS}}_D)$, 
where $D_0$ is an open subset of a space from the class $\mathcal{E}$.\\ 
\end{defi}

\begin{rem} By Remark \ref{manifolds are functored}, for $\mathcal{F}=\bf{LCS}$, every smooth $\mathcal{E}$-manifold is, 
in a natural way, a structure functor-smooth $\mathcal{E}$-manifold. \\
\end{rem}

\begin{thm} Let $\mathcal{E}$ be a class of l.c. models. If $M=(M_0,\mathcal{A})$ and $N=(N_0,\mathcal{B})$ are smooth $\mathcal{E}$-manifolds, and $\Phi:M_0\longrightarrow N_0$ a continuous map, then $\Phi$ is smooth if and only if $\Phi:(M_0,(\mathcal{C}^{\infty})^{\mathbf{LCS}}_M)\longrightarrow (N_0,(\mathcal{C}^{\infty})^{\mathbf{LCS}}_N)$ is a morphism of $\mathbf{LCS}$-functored spaces.
\end{thm}

\noindent {\bf Proof.~} Note that a morphism $(\Phi_0,\Phi^{\sharp})$ between \emph{reduced} functored spaces is completely determined by the underlying continuous map $\Phi_0$, the information contained in the natural transformation $\Phi^{\sharp}$ being that of all possible pullbacks by $\Phi_0$. Consequently, as in the case of \emph{reduced} locally ringed spaces, we write $\Phi$ instead of $(\Phi_0,\Phi^{\sharp})$. Suppose now that $\Phi:M\longrightarrow N$ is smooth. Then, for any open subset $V$ of a locally convex space $F$ and for any open subset $W$ of $N_0$, the pullback $\Phi^*:\mathcal{C}^{\infty}(W,V)\longrightarrow \mathcal{C}^{\infty}(\Phi^{-1}(W),V)$ is well-defined by the chain rule. We claim that $\Phi^{\sharp}:\Phi^{-1}(\mathcal{C}^{\infty})^{\mathbf{LCS}}_N\longrightarrow(\mathcal{C}^{\infty})^{\mathbf{LCS}}_M$, which assigns to every $V$ the sheaf map $\Phi^{\sharp}(V):(\mathcal{C}^{\infty})^{\mathbf{LCS}}_N(V)\longrightarrow\Phi_*((\mathcal{C}^{\infty})^{\mathbf{LCS}}_M(V))$ given by the pullbacks $\Phi^*$, is a natural transformation. Indeed, we have to show that given $V'$ (resp. $V''$) open subset of a locally convex space $F'$ (resp. $F''$), and given a smooth map $\chi:V'\longrightarrow V''$, we have $$(\mathcal{C}^{\infty})^{\mathbf{LCS}}_M(\chi)(\Phi^{-1}(W))\circ \Phi^{\sharp}(V')(W)=\Phi^{\sharp}(V'')(W)\circ (\mathcal{C}^{\infty})^{\mathbf{LCS}}_N(\chi)(W)$$
for every open set $W$ in $N_0$. But this is true since both LHS and RHS, when evaluated at an element $\Psi\in (\mathcal{C}^{\infty})^{\mathbf{LCS}}_N(V')(W)=\mathcal{C}^{\infty}(W,V')$, are equal to $\chi\circ\Psi\circ\Phi$. Thus, $\Phi$ is a morphism of $\mathbf{LCS}$-functored spaces. Now we prove the converse of the theorem. Assume that $\Phi$ induces, via the pullbacks $\Phi^*$, a morphism of $\mathbf{LCS}$-functored spaces. This means that for a germ of a smooth map defined on $N_0$ having values in an open subset of an arbitrary locally convex space, the pullback is the germ of a smooth map on $M_0$. Since this property as well as the condition of smoothness are local, we replace w.l.o.g. $M_0$ and $N_0$ by open subsets $U$ resp. $U'$ in locally convex spaces $E,E'$ in $\mathcal{E}$. Taking now $\Psi=\hbox{id}_{U'}$, the pullback $\Phi^*(\Psi)=\Psi\circ\Phi$ equals $\Phi$ and thus $\Phi$ is smooth. $\hfill\square$\\

The above theorem has the following important consequence.\\

\begin{cor}
Let $\mathcal{E}$ be a class of l.c. models. For every structure functor-smooth $\mathcal{E}$-manifold $(M_0,\mathcal{O}^{\mathbf{LCS}}_M)$, there is a canonical maximal atlas $\mathcal{A}$ on $M_0$ such that $M=(M_0,\mathcal{A})$ is a smooth $\mathcal{E}$-manifold fulfilling $\mathcal{O}^{\mathbf{LCS}}_M=(\mathcal{C}^{\infty})^{\mathbf{LCS}}_M$. Furthermore, the $\mathcal{E}$-manifold $M$ is uniquely determined by the condition $\mathcal{O}^{\mathbf{LCS}}_M=(\mathcal{C}^{\infty})^{\mathbf{LCS}}_M$.
\end{cor}

\noindent {\bf Proof.~} Mutatis mutandis the proof of Corollary \ref{manifolds = structure sheaf-smooth manifolds} shows this corollary as well.  $\hfill\square$\\

One could also take $\mathcal{E}$ to be the class of all complex locally convex spaces. The notion of smooth $\mathcal{E}$-manifold is then replaced by that of {\it complex-analytic locally convex manifold}. To such a manifold $M=(M_0,\mathcal{A})$, one can then associate a {\it complex-analytic locally convex manifold in the structure-sheaf sense} $(M_0,\mathcal{C}^{\omega}_M)$, where $\mathcal{C}^{\omega}_M$ is the sheaf of complex-analytic functions on $M$. One has then the following result.\\

\begin{thm}
Let $M=(M_0,\mathcal{A})$ and $N=(N_0,\mathcal{B})$ be complex-analytic locally convex manifolds, and $\Phi:M_0\longrightarrow N_0$ a continuous map. Then $\Phi$ is complex-analytic if and only if $\Phi:(M_0,\mathcal{C}^{\omega}_M)\longrightarrow (N_0,\mathcal{C}^{\omega}_N)$ is a morphism of locally ringed spaces.\\
\end{thm}

\noindent {\bf Proof.~} By the chain rule, analytic maps are morphisms of locally ringed spaces. Assume now that $\Phi$ is a morphism of locally ringed spaces. Since the condition of analyticity is local, we can proceed as in the proof of Theorem \ref{smooth equals ringed morphism}, which amounts to replace $M_0$ and $N_0$ by open subsets $D_E$ and $D_F$ of complex locally convex spaces $E$ and $F$ respectively. This gives a map $\tilde{\Phi}:D_E\longrightarrow D_F$ which is scalarwise analytic. In particular, viewed as a map from $D_E$ to $F$, $\tilde{\Phi}$ is weakly analytic. By \cite{maz} (Part II, Proposition 1.6), we deduce that $\tilde{\Phi}$ is analytic. It follows that $\Phi$ itself is analytic.   $\hfill\square$\\

\begin{rem}
If, instead of complex-analytic locally convex manifolds, one considers more generally complex-analytic subsets of locally convex spaces
(and analytic spaces modelled on such analytic sets), then pathologies appear. There are examples of reduced analytic sets (in the ringed space-sense) with continuous maps into some complex Banach space which are weakly analytic but not analytic (see e.g. \cite{maz}, pp. 73-80). This phenomenon is avoided by defining analytic sets and spaces as functored spaces in \cite{dou}.\\
\end{rem}

\section{Global characterization of smoothness}

In this section, we apply Theorem \ref{smooth equals ringed morphism} to prove an infinite-dimensional generalization of the classical result which states that a continuous map $\Phi$ between finite-dimensional manifolds $M$ and $N$ is smooth if and only if the pullback by $\Phi$ of every ({\it globally defined}) smooth function on $N$ is smooth (compare, e.g., Lemma 2.2 in \cite{ns}). Our generalization will apply in particular to the case where the target manifold $N$ is the Fr\'echet manifold of smooth maps from a finite-dimensional compact manifold $X$ to a finite-dimensional manifold $Y$. \\

As in the finite-dimensional case, the result depends crucially on the existence of smooth bump functions (and hence is false in the analytic category). As shown in \cite{tho} (but compare also Theorem 16.10 in \cite{km}), smooth bump functions do exist on nuclear Fr\'echet spaces, and this continues to hold true 
for manifolds modeled on nuclear Fr\'echet spaces, provided the manifolds under consideration are regular as topological spaces. More precisely, 
one has the following proposition.\\

\begin{prop}
 Let $N=(N_0,\mathcal{B})$ be a smooth regular manifold modeled on a nuclear Fr\'echet space. For every point $q_0\in N_0$ and for every neighborhood $U$ of $q_0$ in $N_0$, there exists a function $\chi\in\mathcal{C}^{\infty}_N(N_0)$ such that: 
\begin{itemize}
\item $\chi(N_0)\subset [0,1]$, 
\item there exists a neighborhood $V$ of $q_0$ in $U$ such that $\chi_{|V}=1$, 
\item $\chi_{|N_0-U}=0$.
\end{itemize}
\end{prop}

\noindent {\bf Proof.~} cf. \cite{tho}, p. 278. $\hfill\square$\\

With the notations of the preceding section, we have the following result.\\

\begin{thm}\label{local = global}
Let $M=(M_0,\mathcal{A})$ be a smooth $\mathcal{E}$-manifold, and $N=(N_0,\mathcal{B})$ a smooth regular manifold modeled on a nuclear Fr\'echet space $E$ (so $N$ is notably a smooth $\mathcal{M}$-manifold, since $E$ is Mackey-complete). Let $\Phi:M_0\longrightarrow N_0$ be a continuous map. The following statements are equivalent: 
\begin{enumerate}
\item [(i)] $\Phi$ is smooth.
\item [(ii)] For every $f\in\mathcal{C}^{\infty}_N(N_0)$, we have $f\circ\Phi\in\mathcal{C}^{\infty}_M(M_0)$.
\end{enumerate}
\end{thm}

\noindent {\bf Proof.~} Statement (i) immediately implies (ii) by the chain rule, so we only have to show that (ii) implies (i). By Theorem \ref{smooth equals ringed morphism}, this will be true if for any open subset $U$ of $N_0$, and for every $g\in\mathcal{C}^{\infty}_N(U)$, we have $g\circ{\Phi}_{|\Phi^{-1}(U)}\in\mathcal{C}^{\infty}_M(\Phi^{-1}(U))$. Let $p_0$ be an arbitrary point in $\Phi^{-1}(U)$. Let us show that $g\circ\Phi_{|\Phi^{-1}(U)}$ is smooth at $p_0$. By regularity of $N$, there exists a neighborhood $U_0$ of $\Phi(p_0)$ such that $\Phi(p_0)\in U_0\subset\bar{U}_0\subset U$. By the preceding proposition, there exists a function $\chi\in\mathcal{C}^{\infty}_N(N_0)$ such that: 
\begin{itemize}
\item $\chi(N_0)\subset [0,1]$, 
\item there exists a neighborhood $V$ of $\Phi(p_0)$ in $U_0$ such that $\chi_{|V}=1$, 
\item $\chi_{|N_0-U_0}=0$.
\end{itemize}
Define the function $f:N_0\longrightarrow\mathbb{R}$ by $f(q):=\chi(q)\;g(q)$ for every $q\in U$ and $f(q)=0$ for every $q\in N_0-U$. Then $f\in\mathcal{C}^{\infty}_N(N_0)$, and $f_{|V}=g_{|V}$. This implies that $(f\circ\Phi)_{|\Phi^{-1}(V)}=(g\circ\Phi_{|\Phi^{-1}(U)})_{|\Phi^{-1}(V)}$. Now $\Phi^{-1}(V)$ is a neighborhood of $p_0$ (since $\Phi(p_0)\in V$). Smoothness of $f\circ\Phi$ at $p_0$ implies now smoothness of $g\circ\Phi_{|\Phi^{-1}(U)}$ at $p_0$. $\hfill\square$\\

Finally, we claim that the above theorem applies in particular for $N:=\mathcal{C}^{\infty}(X,Y)$, where $X$ and $Y$ are finite-dimensional manifolds, $X$ being compact. Indeed, recall that $\mathcal{C}^{\infty}(X,Y)$ is naturally equipped with the initial topology corresponding to the injection $\displaystyle \iota:\mathcal{C}^{\infty}(X,Y)\longrightarrow \prod_{k=0}^{\infty}\mathcal{C}(\hbox{T}^kX,\hbox{T}^kY)$. (If $\pi_k$ is the canonical projection of the preceding product onto $\mathcal{C}(\hbox{T}^kX,\hbox{T}^kY)$ (which has the compact-open topology), and if\\ $\hbox{T}^k:\mathcal{C}^{\infty}(X,Y)\longrightarrow$ $\mathcal{C}(\hbox{T}^kX,\hbox{T}^kY)$ is the map sending each smooth function $f:X\longrightarrow Y$ to its $k^{th}$-order tangent map $\hbox{T}^kf:\hbox{T}^kX\longrightarrow\hbox{T}^kY$, then $\iota$ is the unique map such that $\pi_k\circ \iota=\hbox{T}^k$ for all $k$.) Since $\hbox{T}^kY$ is metrizable, it is a fortiori regular. It follows that $\mathcal{C}(\hbox{T}^kX,\hbox{T}^kY)$ is also regular (in the compact-open topology). But then $\displaystyle \prod_{k=0}^{\infty}\mathcal{C}(\hbox{T}^kX,\hbox{T}^kY)$ is regular as well, which in turn implies that $\mathcal{C}^{\infty}(X,Y)$ is regular, since $\iota$ is injective. On the other hand, local models for the manifold $\mathcal{C}^{\infty}(X,Y)$ are Fr\'echet spaces of the form $\Gamma_{\mathcal{C}^{\infty}}(X,f^*\hbox{T}Y)$, for $f\in\mathcal{C}^{\infty}(X,Y)$ (compare, e.g. Examples 4.1.2 and 4.1.3 in \cite{ham}). But spaces of smooth sections of vector bundles over a compact manifold are nuclear (cf. the remark below). This establishes our claim.\\

\begin{rem}
Since we could not find a reference for the -certainly folkloristic- fact that the space of smooth sections of a vector bundle $E$ over a compact manifold without boundary 
$X$ is nuclear, we would like to sketch briefly two possible proofs of this fact. One can proceed geometrically by embedding $X$ in a torus $\mathbb{T}^N$ (for $N$ sufficiently large). Then, for $\mathbb{K}=\mathbb{R}$ or $\mathbb{C}$, the space $\mathcal{C}^{\infty}(X,\mathbb{K})$ becomes the quotient of the nuclear space $\mathcal{C}^{\infty}(\mathbb{T}^N,\mathbb{K})$ by a closed subspace, hence is nuclear (see e.g. \cite{tre}). But any vector bundle $E$ over $X$ can be embedded in a trivial vector bundle, and then $\Gamma_{\mathcal{C}^{\infty}}(X,E)$ becomes a subspace of $\mathcal{C}^{\infty}(X,\mathbb{K})\otimes\mathbb{K}^r$ for some $r$, which implies the nuclearity of $\Gamma_{\mathcal{C}^{\infty}}(X,E)$ (compare again \cite{tre}). It could be interesting however to note that a direct analytic approach is possible, in which one generalizes to $\Gamma_{\mathcal{C}^{\infty}}(X,E)$ the standard argument that shows that the space of smooth functions on a torus is nuclear. 
More precisely, choosing a Riemannian metric $g$ on $X$, as well as a bundle metric and a connection on $E$, let $\Delta:\Gamma_{\mathcal{C}^{\infty}}(X,E)\longrightarrow \Gamma_{\mathcal{C}^{\infty}}(X,E)$ be the corresponding Bochner Laplacian. Then $\Gamma_{\mathcal{C}^{\infty}}(X,E)$ is contained in the Hilbert space $\Gamma_{L^2}(X,E)$, and $\langle\Delta^k \psi,\psi\rangle_{L^2}<\infty$ for all $\psi \in \Gamma_{\mathcal{C}^{\infty}}(X,E)$ and for all $k\in\mathbb{N}$. If $(\lambda_n)_{n\in\mathbb{N}}$ is the increasing sequence of eigenvalues of $\Delta$, and $(\psi_n)_{n\in\mathbb{N}}$ an orthonormal basis of $\Gamma_{L^2}(X,E)$ made of smooth eigensections of $\Delta$, this gives: $$\sum_{n=0}^{\infty}(\lambda_n)^k\; |\langle \psi,\psi_n\rangle_{L^2}|^2<\infty\quad\quad\forall k\in\mathbb{N}\, .$$
Now given $\alpha\in\mathbb{N}$, choose $k\in\mathbb{N}$ such that $k>\frac{\alpha d}{2}$, where 
$d:=\hbox{dim}_{\mathbb{R}}X$. Weyl's asymptotic formula for generalized Laplacians acting on sections of vector bundles, as given, e.g.,  in Corollary 2.43 of \cite{bgv}, is equivalent to
$$\lambda_n\;\sim\;\frac{4\pi^2}{[\hbox{rank}(E)\;\hbox{vol}(\hbox{\upshape{B}}_{\mathbb{R}^d})\;\hbox{vol}_g(X)]^{\frac{2}{d}}}\;n^{\frac{2}{d}}\quad\quad\hbox{as }\;n\to\infty \, .$$
Accordingly, we have $(\lambda_n)^k\sim C n^{\frac{2k}{d}}$ for some constant $C$. It follows easily that $n^{\alpha}\ll(\lambda_n)^k$ as $n\to\infty$, and so: 
$$\sum_{n=0}^{\infty}n^{\alpha}\; |\langle \psi,\psi_n\rangle_{L^2}|^2<\infty\quad\quad\forall\alpha\in\mathbb{N}  \, .$$ 
Thus, the generalized Fourier transform $\psi\mapsto (\langle \psi,\psi_n\rangle_{L^2})_{n\in\mathbb{N}}$ defines an injective continuous linear map from $\Gamma_{\mathcal{C}^{\infty}}(X,E)$ to the space 
$s(\mathbb{N})$ of rapidly decreasing sequences, which is the prototype of nuclear spaces. This map is easily seen to be surjective as well, and by the open mapping theorem for Fr\'echet spaces, it follows that it is a linear homeomorphism. Using Theorem 51.5 of \cite{tre} this implies the nuclearity of $\Gamma_{\mathcal{C}^{\infty}}(X,E)$.
\end{rem}

\vskip2cm

\end{document}